\documentclass[11pt,twoside]{article}
\usepackage{graphics,latexsym,amssymb,amsmath}
\RequirePackage{graphics,epsfig,psfrag,picins} 
\def\RR{{\rm I\!R}} 
\def\pn{\medskip\par\noindent}
\def\Frac#1#2{{\displaystyle{{#1}\over{#2}}}}

\def\[#1\]{\begin{eqnarray*}#1\end{eqnarray*}}
\def\$#1\${\begin{eqnarray}#1\end{eqnarray}}

\def\phi{\varphi}

\def\eps{\varepsilon}
\def\bref#1{{\bf \ref{#1}}}
\newcommand{\Pf}{{\em Proof}. }
\newcommand{\prbox}{\EPf}
\newcommand{\EPf}{\hbox{}\hfill$\Box$\vspace{.5cm}}

\date{\today}

\def\Frac#1#2{{\displaystyle{\frac{#1}{#2}}}}
\def\phi{\varphi}

\newtheorem{definition}{Definition}
\newtheorem{remark}{Remark}
\newtheorem{lemma}{Lemma}
\newtheorem{thm*}{Theorem}
\date{\today}

\begin{document}
\pagestyle{myheadings}
\markboth{P. -V. Koseleff, D. Pecker, {\em On polynomial Torus
    Knots}}{Accepted for publication in J. of Knot Th. and Ram.}

\title{On polynomial Torus Knots}

\author{P. -V. Koseleff, D. Pecker\medskip\\
Universit{\'e} Pierre et Marie Curie\\
4, place Jussieu, F-75252 Paris Cedex 05 \\
e-mail: {\tt\{koseleff,pecker\}@math.jussieu.fr}}

\maketitle

\begin{abstract}
We show that no torus knot of type $(2,n)$, $n>3$ odd, can be obtained from a
polynomial embedding $t \mapsto ( f(t), g(t), h(t) )$ where
$(\deg(f),\deg(g))\leq (3,n+1) $. Eventually, we give explicit examples
with minimal lexicographic degree.
\end{abstract}
{\bf keywords:} {Knot theory, polynomial curves, torus knots, 
parametrized space curve} \\
{\bf Mathematics Subject Classification 2000:} 14H50, 12D10, 26C10, 57M25
\section{Introduction}
The study of non compact knots began with Vassiliev \cite{Va}.
He proved that any non-compact knot type can be obtained from
a polynomial embedding
$ \  t \mapsto ( f(t), g(t), h(t) ),  t \in \RR .$
The proof uses
Weierstrass approximation theorem on a compact interval,
the degrees of the polynomials may be quite large,
and the plane projections of the polynomial knots quite complicated.

Independently, Shastri \cite{Sh} gave a detailed proof of this theorem,
he also gave simple polynomial parametrizations of the trefoil and of the
figure eight knot.

This is what motivated  A. Ranjan and Rama Shukla \cite{RS}
 to find  small degree parametrizations of the simplest knots,
 the torus knots of type $(2,n), n$ odd, denoted by $K_n$.
They proved that these knots can be attained from polynomials of degrees
 $ (3, 2n-2, 2n-1).$
In particular,
they obtain a parametrization of the trefoil $K_3$ analogous to Shastri's one.
They also asked the natural question which is to find the minimal
degrees of the polynomials representing a general torus knot of a
given type (there is an analogous question in Vassiliev's paper
\cite{Va}).
 
The number of  crossings of a plane projection of $K_n$  is at least
$n$ (Bankwitz theorem, see \cite{Re}). 
It is not difficult to see, using B{\'e}zout theorem, that this plane
curve cannot be parametrized by  polynomials of degrees smaller than $(3,n+1).$

Naturally, Rama Mishra (\cite{Mi}) asked whether it was possible  
to parametrize the knot $ K_{n } $ by polynomials of degrees  $ (3,
n+1, m)$ when $ n \equiv 1$, or $0 \mod 3.$ 

In this paper, we shall prove the following result
\pn
{\bf Theorem}.
{\it 
If $n \ne 3 $ is odd, the torus knot  $ K_{n} $ cannot be represented by
polynomials of degrees $(3, n+1, m).$}
\pn
Our method is based on the fact that 
all plane projections of  $K_n$ with the minimal number
$n$ of crossings have essentially the same diagram.
This is a consequence of the now solved classical Tait's
conjectures \cite{Mu,Ka,Pr,MT}.
This allows us to transform our problem into a problem of real
polynomial algebra.

As a conclusion, we give explicit parametrizations
of    $K_3,  K _{ 5}  $ and $ K _{ 7 }.$ By our result, they are of
minimal degrees. We also give an explicit parametrization of $K_9$ 
with a plane projection possessing the minimal number of crossing
points. This embedding is of smaller degree than those already known. 
\section{The principal result}
If $n$ is odd, the torus knot $K_n$ of type $(2,n)$ is the
boundary of a  Moebius band  twisted $n$ times (see \cite{Re,Ka,St}).
\begin{figure}[th]
\centerline{
\psfig{file=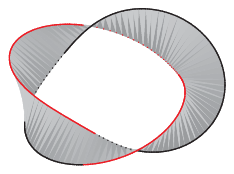} \quad
\psfig{file=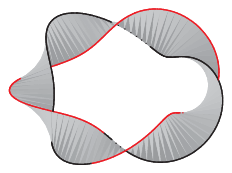} \quad
\psfig{file=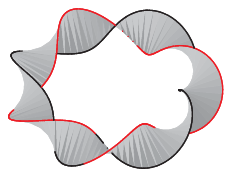}}
\vspace{10pt}
\centerline{
\psfig{file=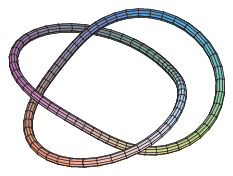} \quad
\psfig{file=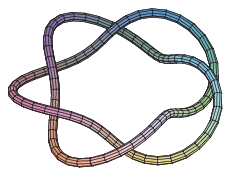} \quad
\psfig{file=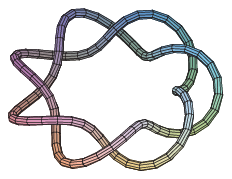}}
\vspace*{8pt}
\caption{$K_n$, $n=3,5,7$.\label{fig1}}
\end{figure}
The recently proved Tait's conjectures allow us to characterize 
plane projections of $K_n$ with the minimal number of crossings.
\begin{lemma}\label{st}
Let ${\cal C}$ be a plane curve with $n$ crossings parametrized by
$ {\cal C }(t) = ( x(t), y(t)) .$  If $ {\cal C} $ is the projection of
a knot $K_n$ then there exist real numbers
$ s_1 < \cdots < s_n < t_1 < \cdots < t_n ,$ such that
${\cal C }(s_i)={\cal C}(t_i).$
\end{lemma}
\Pf
Let  ${\cal C}$ be a plane projection of a knot of type $K_n$
with the minimal number $n$ of crossings. 

Using the Murasugi's theorem B (\cite{Mu}) which says that a minimal
projection of a prime alternating knot is alternating,
we see that ${\cal C } $ is  alternating.

Then the Tait's flyping conjecture, proved by Menasco and Thistlethwaite
(\cite{MT,Pr}), asserts that ${\cal C}$ is related to the standard 
diagram of $K_n$ by a sequence of flypes.
Let us recall that a flype is a transformation most clearly described by
the following picture.
\begin{figure}[th]
\centerline{
{\scalebox{.8}{\includegraphics[width=10cm]{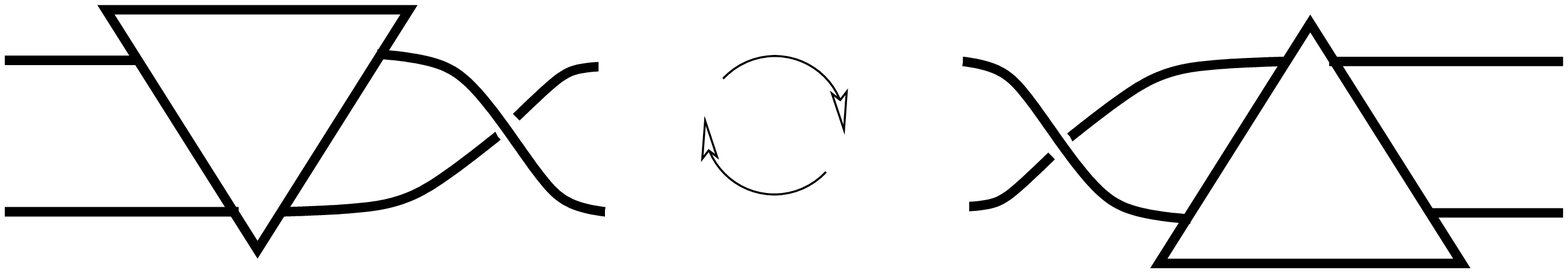}}}}
\caption{A flype\label{fig2}}
\end{figure}
\pn
The standard diagram ${\cal S}_0$ of $K_n$ has the property (cf \cite{Re}) 
that there exist real numbers $ s_1< \cdots <s_n < t_1 < \cdots < t_n $ 
such that $ {\cal S}_0 (s_i) = {\cal S}_0 (t_i) .$ It is alternating. 

Let ${\cal S}$ be a diagram with real parameters $ s_1< \cdots <s_n <
t_1 < \cdots < t_n $  such that $ {\cal S} (s_i) = {\cal S } (t_i)$,
and let us perform a flype of a part $B$ of ${\cal S}$
\begin{figure}[th]
\psfrag{A}{}
\psfrag{B}{}
\psfrag{C}{}
\psfrag{D}{}
\psfrag{E}{}
\psfrag{F}{}
\psfrag{g}{$s$}
\psfrag{h}{$t$}
\psfrag{U}{{A}}
\psfrag{V}{{B}}
\psfrag{W}{{C}}
\centerline{
{\scalebox{.8}{\includegraphics[width=8cm]{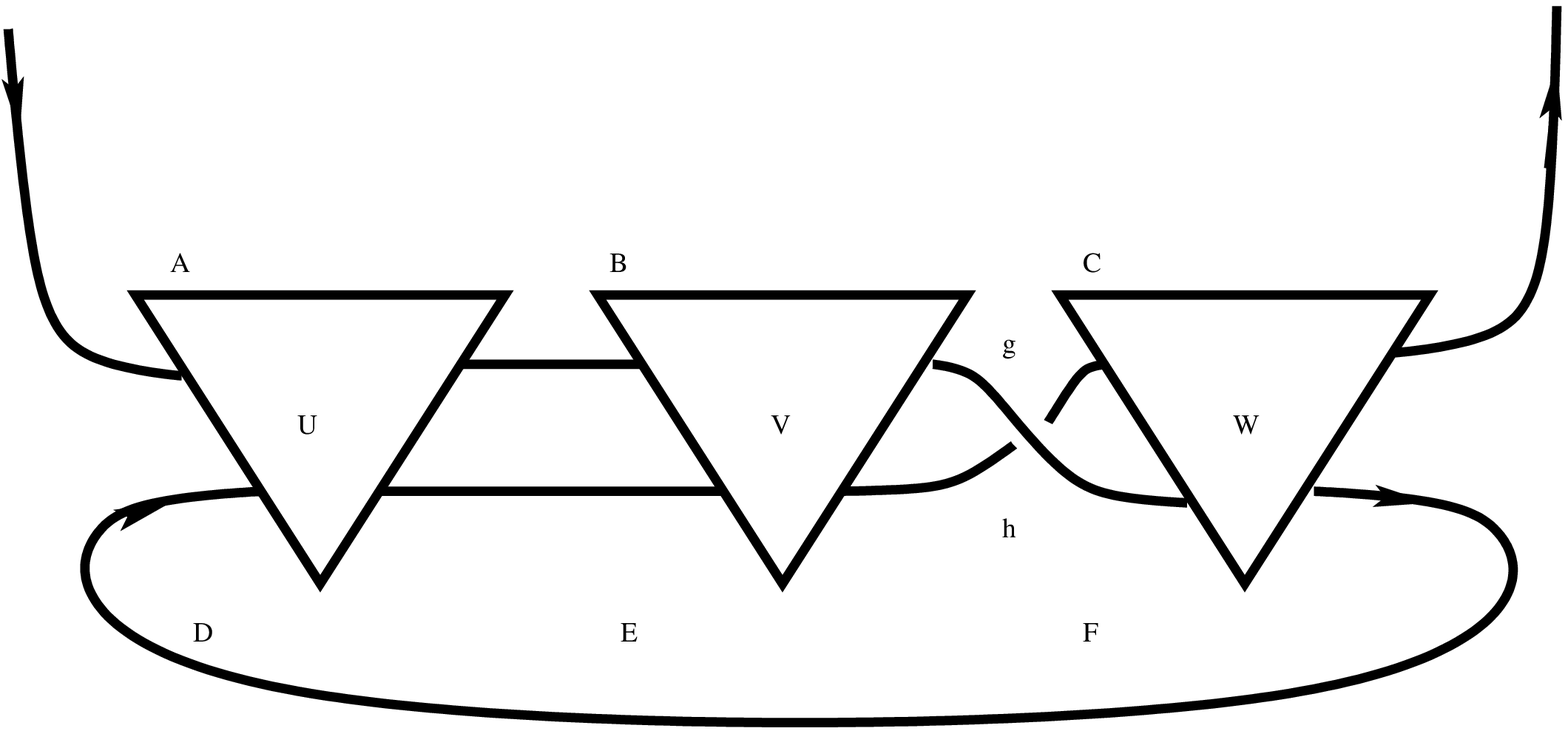}}}
\quad
\psfrag{g}{$s'$}
\psfrag{h}{$t'$}
{\scalebox{.8}{\includegraphics[width=8cm]{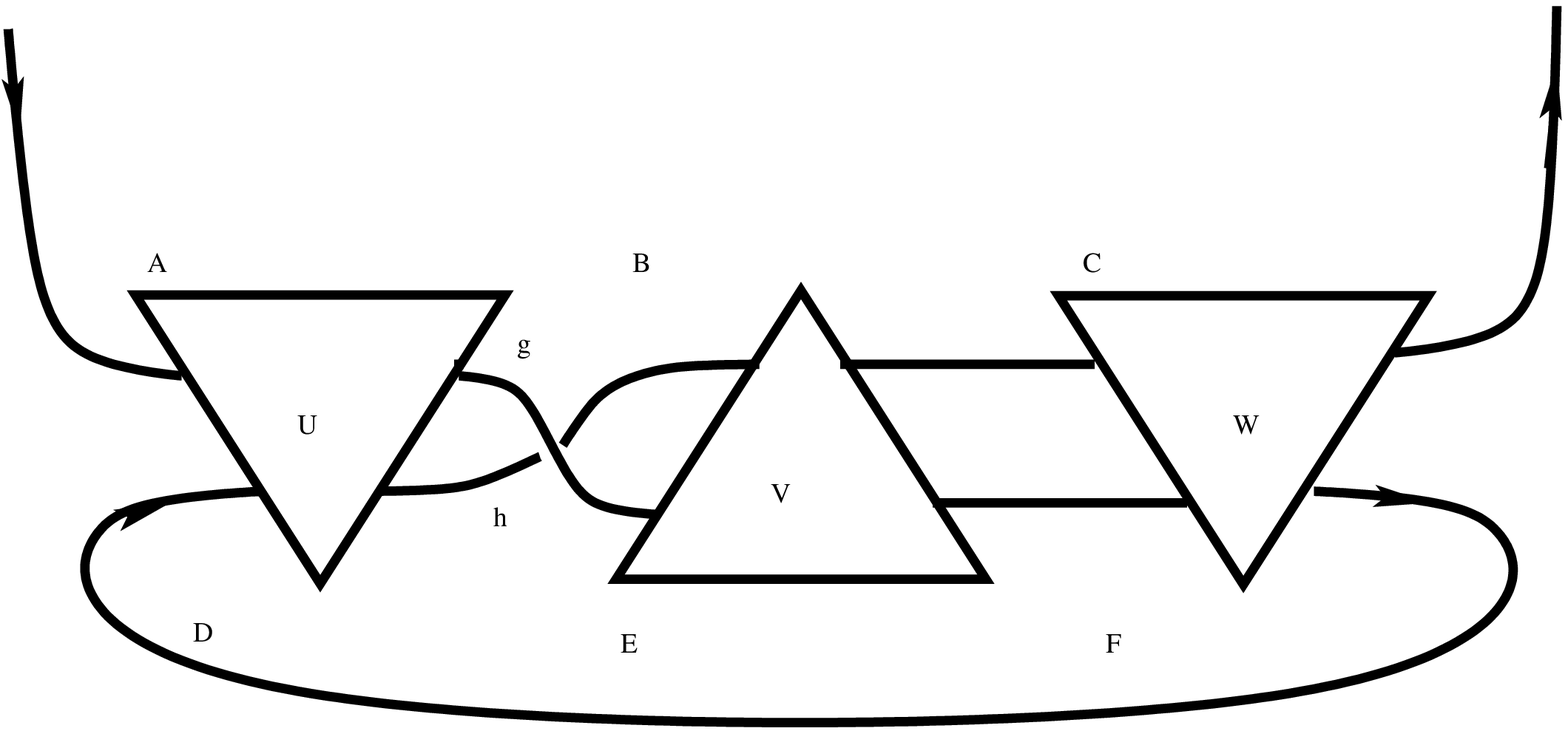}}}}
\caption{Flype on the part $B$\label{fig3}}
\end{figure}
\par\noindent
For any $(a,b,c) \in A \times B \times C$ we have 
$$
s_a < s_b < s < s_c < t_a < t_b < t < t_c.
$$
After the flype on $B$, we have new parameters corresponding to
the crossing points satisfying 
$$
s_a' < s' < s_b' <  s_c' < t_a' <t'<  t_b' < t_c'.
$$
The transformed diagram ${\cal S}'$ has the same property: there exist
real parameters $ s_1 < \cdots < s_n < t_1 < \cdots < t_n ,$ such that
${\cal S}'(s_i)={\cal S}'(t_i).$

So then, after any sequence of flypes, the transformed diagram will
have the same property.
\EPf

In this paper we shall consider polynomial knots, that is to say, polynomial
embeddings
$ {\RR} \longrightarrow {\RR}^3, \   t \mapsto ( x(t), y(t), z(t)). $
Polynomial knots are non-compact subsets of $\RR^3$.
The closure of a polynomial knot in the one point compactification
$ {\bf S}^3$ of the space $ {\RR}^3$ is an ordinary knot
(see \cite{Va,Sh,RS} and figures at the end).

\begin{lemma}
Let ${\cal C}$ be a plane polynomial curve
with $n$ crossings parametrized by $ {\cal C} (t) = (x(t),y(t)).$
Suppose that ${\cal C}$ is the projection of $K_n$ and
$ \deg x(t) \le \deg y(t) .$
Then we have $ \deg x(t) \ge 3 .$  If $ \deg x(t) =3 ,$
then $ \deg  y(t) \ge n+1 .$
\end{lemma}
\Pf
$x(t)$ must be non-monotonic, so $ \deg x(t) \ge 2 .$
Suppose that $x(t)$ is of degree~2. Then $x(t_i) = x(s_i)$
implies that $ t_i+s_i$ is constant, and so the parameter values
corresponding to the crossing points are ordered as
$$ s_1< \cdots <s_n < t_n < \cdots < t_1. $$
We have a contradiction according to lemma \bref{st}
\pn
Suppose now that $ \deg x(t) =3$. The crossing points of the curve 
${\cal C}$ correspond to parameters $(s,t)$, $s\not = t$, that are
common points of the curves of degrees $2$ and $\deg y(t) -1$:
$$
\Frac{x(t)-x(s)}{t-s}= 0, \; 
\Frac{y(t)-y(s)}{t-s}= 0.
$$
By B{\'e}zout theorem (\cite{Fi}), the number of such points are at most
$2 \times (\deg y(t) -1)$.
$(s,t)$ and $(t,s)$ are distinct points and correspond to the
same crossing point. So, the curve $ {\cal C}$ has at most 
$\deg y(t) -1$ crossing points, and this implies that 
$ \deg y(t) \ge n+1 .$
\EPf

\section{Proof of the main result}
Our proof makes use of Chebyshev (monic) polynomials.
\subsection{Chebyshev Polynomials}
\begin{definition}
If $t= 2 \cos \theta$, let $T_n(t)= 2 \cos ( n \theta )$ and
$V_n (t) = \Frac{\sin ((n+1) \theta)}{\sin \theta }$.
\end{definition}
\begin{remark} $T_n$ and $V_n$ are both monic and have degree $n$.
We have 
\begin{equation}
V_0=1,\quad V_1 = t , \quad V_{n+1} = t\, V_n  - V_{n-1}.\label{recv}
\end{equation}
We have also 
\begin{equation}
T_0 = 2,\quad T_1 = t, \quad 
T_{n+1} = t \, T_n - T_{n-1}. \nonumber
\end{equation}
For $n \geq 2$, let $V_n = t^n + a_n t^{n-2} + b_n t^{n-4} + \cdots$.
Using recurrence formula 
\bref{recv}, we get 
\begin{equation*}
a_{n+1} = a_n - 1, \, b_{n+1} = b_{n} - a_{n-1}
\end{equation*}
so by induction, 
\begin{equation}
V_n = t^n - (n-1) t^{n-2} + \Frac 12 (n-2)(n-3) t^{n-4} + \cdots. 
\label{dev}
\end{equation}
\end{remark}
We shall also need the following lemmas which will be proved in the
next paragraph.
\pn
{\bf Lemma A.}
{\it 
Let $s \not = t$ be real numbers such that  $T_3(s)= T_3(t).$ For
any integer $k$, we have
\begin{equation}
\Frac{T_k (t) - T_k(s)}{t - s} =\Frac{2}{\sqrt 3} \sin \Frac{k
\pi}3 \, V_{k-1}(s+t). \nonumber
\end{equation}
}
\pn
{\bf Lemma B.}
{\it Let $n\ge 3 $ be an integer.\\
Let $s_1<s_2 < \cdots < s_n$ and $t_1 < \cdots < t_n$ be real
numbers such that  $T_3(s_i)= T_3(t_i).$ Let $u_i = t_i + s_i$. We
have
\begin{equation*}
\sum_{i=1}^n u_i^2 \leq  n + 4 , \ \sum_{i=1}^n u_i^4 \leq  n + 22. 
\end{equation*}
}
\subsection{Proof of the theorem}
\Pf
We shall prove this result by reducing it to a contradiction. Suppose
the plane curve $ {\cal C}$ parametrized by $ x= P(t), \  y =
Q(t)$ where $ \deg P =3, \  \deg Q = n+1$ is a plane projection
of $K_{n}.$ 

By translation on $t$, one can suppose that $P(t)=t^3-\alpha t + \beta$.
If the polynomial $P$ was monotonic, $ {\cal C}$ would
have no crossings, which is absurd. Therefore $\alpha>0$.
Dividing $t$ by $\rho=\sqrt 3/\sqrt{\alpha}$, one has $P(t) = \rho^3 (t^3
- 3t) + \mu.$ By translating the origin and scaling $x$, one can now
suppose that $P(t) = t^3 - 3t = T_3(t)$. 

By translating the origin and scaling $y$,  
we can also suppose that $Q(t)$ is monic and write
\begin{equation}
P(t)= T_3(t), \quad 
Q(t) = T_{n+1}(t) + a_n T_n(t) + \cdots + a_1 T_1(t).
\nonumber
\end{equation}
By B{\'e}zout
theorem, the curve  ${\cal C}$ has at most  $(3-1)(n+1-1)/2 = n $
double points. As it has at least $n$ crossings, we see that
it has exactly  $n$ crossings and therefore is a minimal crossing
diagram of $ K_{n}$.
According to the lemma {\bf \bref{st}}, there  exist   real numbers
$ s_1 < \cdots <s_n , \  t_1 < \cdots < t_n, \  s_i < t_i $
such that $P(s_i)= P(t_i), \   Q(s_i)= Q (t_i) .$
Let $u_i = t_i+s_i, \ 1 \leq i \leq n$. We have
\begin{equation}
\Frac{Q(t_i) - Q(s_i)}{t_i - s_i} =
\Frac{T_{n+1} (t_i) - T_{n+1} (s_i)}{t_i -s_i} + \sum_{k=1}^n a_k
\Frac{T_k(t_i) - T_k(s_i)}{t_i - s_i}. \nonumber
\end{equation}
so by lemma {\bf A},
$u_1, \ldots, u_n$ are the distinct roots of the polynomial
\begin{equation}
R(u) = \eps_{n+1}
\, V_{n}(u)
+
\sum_{k=1}^{n} a_k \eps_k \,
V_{k-1}(u),
\end{equation}
where 
$$\eps_k = 
\Frac{2}{\sqrt 3} \sin \Frac{k\pi}3.$$ 
\begin{remark}
Note that $\eps_k = V_{k-1}(1)$ is the
6-period sequence $\eps_0=0$, $\eps_1=1$, $\eps_2=1$, $0$, $-1$, $-1, \ldots$.
\end{remark}
We have to consider several cases.
\pn 
$\rhd$ {\bf Case $\mathbf{n \equiv 2  \mod 3}$}. $\eps_{n+1} = 0$ and 
$R(u)$ has degree at most $n-1$. This is a contradiction.  
\pn
$\rhd$
{\bf Case ${\mathbf{n \equiv 1  \mod 3}}$}.
In this case, $n \equiv 1 \mod 6$ and $\eps_{n+1} = \eps_{n}=1$,
$\eps_{n-1} = 0$. Thus $R(u)$ can be written as
\begin{eqnarray*}
R(u)&=&V_n(u) + a_{n} V_{n-1} (u) - a_{n-2} V_{n-3}(u) - \cdots + 
a_2 V_1(u) + a_1 \nonumber \\
&=& u^{n} + a_{n} u^{n-1} - (n-1) u^{n-2} + \cdots .
\end{eqnarray*}
using equation {\bf \bref{dev}}.
Therefore we get
\begin{equation}
\sum_{1 \leq i \leq n} u_i = - a_{n}, \
\sum_{1 \leq i < j \leq n} u_i u_j = - ( n-1) ,\nonumber
\end{equation}
and then
\begin{equation}
\sum _{i=1} ^n u_i^2 = \left ( \sum _{i=1}^n u_i \right ) ^2 - 2 \sum
_{1 \leq i < j \leq n} u_i u_j = 
a_{n}^2 + 2 ( n-1)  \ge 2 ( n-1).\nonumber
\end{equation}
According to lemma {\bf B} we also have $\sum_{i=1}^n u_i^2 \le n+4,$ we
get a contradiction for $n>6.$
\pn
$\rhd$
{\bf Case $\mathbf{ n \equiv 0 \mod 3}$.}
In this last case we have $n=3 \mod 6$, so $\eps_{n+1}=-1$, $\eps_n=0$
and $\eps_{n-1}=1$, so  
\begin{eqnarray*}
- R(u)&=& V_n(u)- a_{n-1} V_{n-2} (u) - a _{n-2} V_{n-3} (u) + \cdots
-a_2 V_1(u)- a_1.
\end{eqnarray*}
Let $\sigma_i$ be the  coefficients of  
$$- R(u) = u^n + \sum_{k=1}^n (-1)^k \sigma_k u^{n-k} .$$
From the equation \bref{dev}, we see that
\begin{equation*}
\sigma_1 = 0, \;
\sigma_2= - (a_{n-1} +n -1) , \;
\sigma_4 = (n-3)\,a_{n-1}+{(n-2)(n-3) \over 2}.
\end{equation*}
Let $S_k$ be the  Newton sums  $\sum _{i=1}^n u_i^k$ of the roots of
the polynomial $R$. Using the classical Newton formulas (\cite{FS}), we obtain
\begin{equation*}
S_1=\sigma_1 = 0, \  S_2= \sigma_1^2 - 2 \sigma_2 = -2 \sigma_2, \  S_4
= 2 \sigma_2^2 - 4 \sigma_4, 
\end{equation*}
and then
\begin{equation*}
S_4= 2 ( a_{n-1} +2)^2 + 6 n - 18 \geq 6n -18.
\end{equation*}
By the lemma {\bf B}, we deduce that
$ 22+n \geq 6 n - 18, $ i.e. $n \le 8$ so $n=3$.
\EPf

\subsection{Proof of lemmas A and B}
We shall use the following lemma
\begin{lemma}[Lissajous ellipse]\label{ellipse}
Let $s \ne t$ be complex numbers such that $$ T_3(t)=T_3(s) .$$
There exists a complex number   $ \alpha $ such that
$$s= 2 \cos ( \alpha + \pi/3 ), \
t= 2\cos ( \alpha - \pi/3 ) .$$
Furthermore, $\alpha$ is real if and only if $s$ and  $t$ are both real,
and then $ t>s $ if and only if $\sin \alpha >0 .$
\end{lemma}
\Pf
We have
\begin{equation} {T_3(t)- T_3(s) \over t-s } = t^2 +s^2 +st -3. \label{ell}
\end{equation}
Then, if $T_3(t)=T_3(s),$ $ t\ne s ,$ we get
\begin{equation*} 
{3 \over 2 } (t+s)^2 + {1 \over 2} (t-s)^2= 2( t^2+s^2+st) = 6. 
\end{equation*}
\parpic(0cm,7cm)(0cm,6cm)[lt]{
\psfrag{u}{$(s,t)$}
\psfrag{v}{$(s,2\cos \alpha)$}
{\scalebox{.8}{\includegraphics[width=7cm]{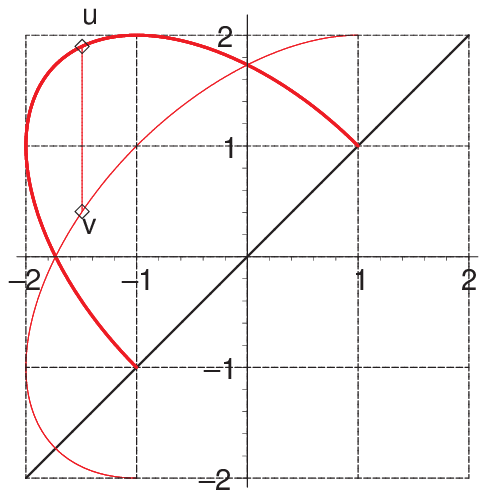}}}}
\noindent
That means
\begin{equation*} \left [ \Frac{t+s}2 \right ]^2 +
 \left [\Frac{t-s}{2 \sqrt 3} \right ] ^2 = 1 .
\end{equation*}
Then there exists a complex number $ \alpha$ such that
\begin{equation*}
\cos\, \alpha = \Frac{t+s}2, \   \sin\, \alpha = \Frac{t-s}{2 \sqrt 3},
\end{equation*}
that is
\begin{equation*}
t = 2 \cos ( \alpha - \pi/3), \
s = 2 \cos ( \alpha + \pi/3) .
\end{equation*}
$ \alpha$ is real if and only if $\cos\, \alpha$ and $\sin\, \alpha$
are both real that is to say, iff $ s$ and  $t$ are real. In this case:
$t>s \Leftrightarrow \sin \alpha >0.$
\EPf
\pn
In order to prove lemma {\bf B}, we shall use the following lemma which
describes the geometrical configuration. 
Let us denote
$s(\alpha)=2 \cos ( \alpha + \pi/3)$ and $
t(\alpha) =2 \cos ( \alpha - \pi/3) .$
\begin{lemma}\label{ineq}
Let $\alpha, \alpha'\in  [0,  \pi ] $  be such that
$ s( \alpha) < s(\alpha'),$ and
$ t( \alpha) < t ( \alpha ') .$
Then $ \alpha > \alpha'$
and $  \Frac{2\pi}{3} > \Frac{\alpha + \alpha'}2 > \Frac{\pi}3 .$
\end{lemma}
\Pf
We have $2 \cos \alpha = s(\alpha)+t(\alpha) < s(\alpha')+t(\alpha') =
2 \cos \alpha'$ so $\alpha > \alpha'.$
\begin{equation*}
\begin{array}{rcccl}
t( \alpha') - t ( \alpha)&=&
4 \sin ( \Frac{ \alpha + \alpha'}2-\Frac{\pi}3 )\cdot
\sin (\Frac{\alpha - \alpha'}2) >0, \\
s( \alpha ') - s( \alpha)& =&
4 \sin ( \Frac{ \alpha + \alpha'}2+\Frac{\pi}3)\cdot
\sin ( \Frac{\alpha - \alpha'}2)  >0.
\end{array}
\end{equation*}
From $ \alpha>\alpha '$ we get 
$0 < \Frac{\alpha + \alpha'}2 - \Frac{\pi}3$ and 
$\Frac{\alpha + \alpha'}2 + \Frac{\pi}3 < \pi$, that is to say 
$$
\Frac{\pi}{3} < \Frac{\alpha + \alpha'}2 < \Frac{2\pi}3 .
$$
\EPf
\pn
{\bf Proof of lemma B}.\\
Let $ s_1 < \cdots <  s_n$ and $ t_1 < \cdots < t_n $ be such that $T_3(s_i)=
T_3( t_i)$. Using lemmas {\bf \bref{ellipse}} and {\bf \bref{ineq}} there are
$ \alpha_1 > \cdots > \alpha_n \in \; ]0, \pi [ $ such that
$ s_i= s( \alpha _i ) , \  t_i = t( \alpha _i )$ and we have
\begin{equation*}
\Frac{2\pi}{3}> \Frac{\alpha_1 + \alpha_2}2 > \alpha_2 > \cdots >
\alpha_{n-1}>\Frac{\alpha_{n-1} + \alpha_n}2 >  \Frac{\pi}3.
\end{equation*}
At least two of the $\alpha_i$'s lie in the intervals $]0,\pi/2]$ or
$[\pi/2,\pi[$. We have only two cases to consider:
$\pi >\alpha_1>\alpha_2\geq \Frac{\pi}{2}$,
or
$\Frac{\pi}{2}\geq \alpha_{n-1} >\alpha_n>0$. 

On the other hand, we get the equality
\begin{equation}
\cos^2 x + \cos^2 y = 1 - \cos^2(x+y) + 2
\cos x \cos y \cos(x+y).
\label{cos}
\end{equation}
\pn
$\rhd$ {\bf Case 1. $\Frac{\pi}{2}\geq \alpha_{n-1} >\alpha_n>0$.}\\
We get $\cos \alpha_n \geq 0$, $ \cos \alpha_{n-1} \geq 0$ and
$\cos(\alpha_{n-1}+\alpha_{n})<-\Frac 12$ so eq. {\bf \bref{cos}} becomes
\begin{equation*}
\cos^2 \alpha_{n-1} + \cos^2 \alpha_n \leq
1 - \cos^2(\alpha_{n-1}+\alpha_{n}) \leq \Frac 34
\end{equation*}
and
\begin{eqnarray*}
\sum _{i=1}^n \cos^2  \alpha_i
&=&
\cos^2 \alpha_1
+
\sum_{i=2}^{n-2} \cos^2 \alpha_i
+
(\cos^2 \alpha_{n-1} + \cos^2 \alpha_n)\nonumber\\
&\leq& 1
+
(n-3) \cdot \Frac 14
+ \Frac 34 = \Frac 14 (n+4),
\end{eqnarray*}
that is
\begin{eqnarray*}
S_2= \sum _{i=1}^n u_i^2 = 
\sum _{i=1}^n \left (2 \cos \alpha _i \right )^2  
\leq n+4.
\end{eqnarray*}
\pn
$\rhd$ {\bf Case 2. $\pi >\alpha_1>\alpha_2\geq \Frac{\pi}{2}$.}\\
We get 
$\cos \alpha_1 \leq 0$, $ \cos \alpha_{2} \leq 0$ and
$\cos (\alpha_{1}+\alpha_{2})<-\Frac 12$ so eq. {\bf \bref{cos}} becomes
\begin{equation*}
\cos^2 \alpha_{1} + \cos^2 \alpha_2 \leq
1 - \cos^2 (\alpha_{1}+\alpha_{2}) \leq \Frac 34
\end{equation*}
and similarly, we get
\begin{equation*}
S_2= \sum _{i=1}^n \left (2 \cos \alpha _i \right )^2 \leq n+4.
\end{equation*}
Analogously, we get
$\cos^4 x + \cos^4 y \leq (\cos^2 x+\cos^2 y)^2$, and we deduce:
\begin{equation*}
S_4 = \sum _{i=1}^n \left ( 2 \cos \alpha_i \right )^4 \leq n+22.
\end{equation*}
\prbox
\pn
{\bf Proof of lemma A}.\\
Let $s<t $ be real numbers such that $T_3(s)=T_3(t).$ 
According to the ellipse lemma \bref{ellipse},
there exists a real number $ \alpha$ such that 
$$t = 2 \cos ( \alpha - \pi/ 3 ), \  
s= 2 \cos (\alpha + \pi/3) .$$
We have
$ s+t = 2 \cos \alpha$ and $t-s = 4 \sin \Frac{\pi}3 \sin \alpha$, so
\begin{eqnarray*}
\Frac{T_k(t)-T_k(s)}{t-s} &=& \Frac{
 2 \bigl( \cos k ( \alpha - \pi/3) -
\cos  k ( \alpha + \pi/3 ) \bigr)}{4 \sin \Frac{\pi}3 \sin \alpha}
=
\Frac{\sin k \alpha \cdot {\sin \Frac{k\pi}{3}}}
{\sin \alpha \cdot \sin\Frac{\pi}{3}}
\\
&=&
\Frac{2}{\sqrt 3} \sin\Frac{k\pi}{3}\, V_{k-1}(2\cos \alpha).
\end{eqnarray*}
\prbox
\section{Parametrized models of  $K_3, \; K_{5}, \; 
K_{7}$ and $K_9$}
We get parametrizations of $K_n$: 
${\cal C} = (x(t), y(t), z(t)),$
with $n$ crossings obtained for parameter values satisfying the
hypothesis of lemma \bref{st}. According to lemma {\bf A}, we choose 
$n$ distinct points $-1 \leq u_1 <\cdots
<u_n \leq  1$. 
We look for $Q_1$ and $Q_2$ of minimal degrees, such that 
$$
R_1(t+s) = \Frac{Q_1(t)-Q_1(s)}{t-s}, \  
R_2(t+s) = \Frac{Q_2(t)-Q_2(s)}{t-s} \  
$$
satisfy, for $i = 1, \ldots, n, $
$$
R_1(u_i) = 0, \, R_2(u_i) = (-1)^i.
$$
We then choose $y(t)=Q_1(t)$ and $z(t)=Q_2(t)$. 
We also add 
some linear combinations of $T_{6i}$ efficiently.
We then obtain a knot whose projection is alternating, when $R_1$ has
no more roots in $[-2,2]$.  As we have chosen symmetric $u_i$'s, all
of our curves are symmetric with respect to the $y$-axis. 

\subsection{Parametrization of $K_3$}
We can parametrize  $K_3$ by $x=T_3(t), y= T_4(t),
z= T_5(t).$ It is a Lissajous space curve (compare \cite{Sh}).
\begin{figure}[th]
\begin{center}
\begin{tabular}{ccc}
\psfig{file=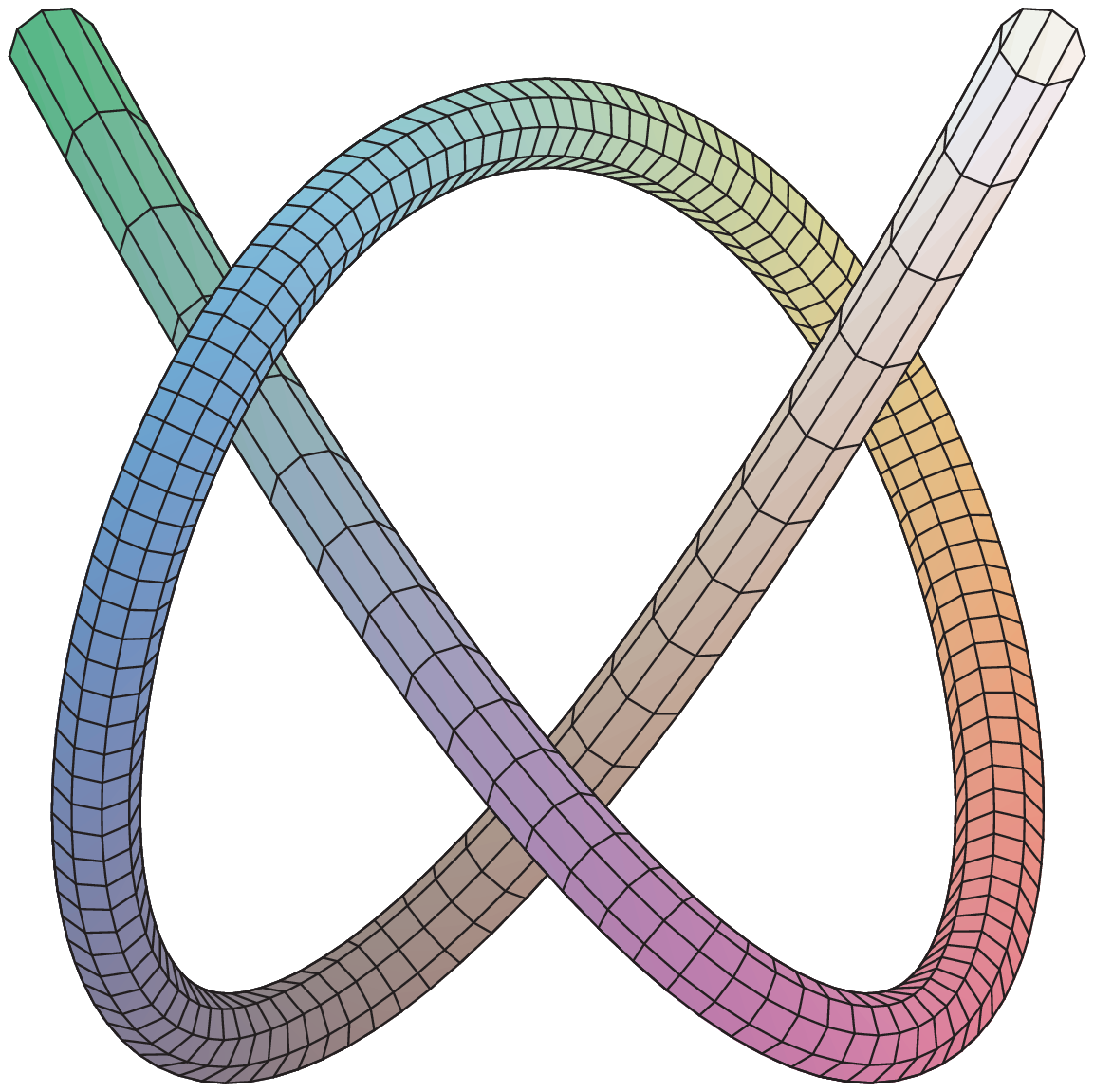,width=4cm}&
\hbox{}\hspace{2cm}\hbox{}&
\psfig{file=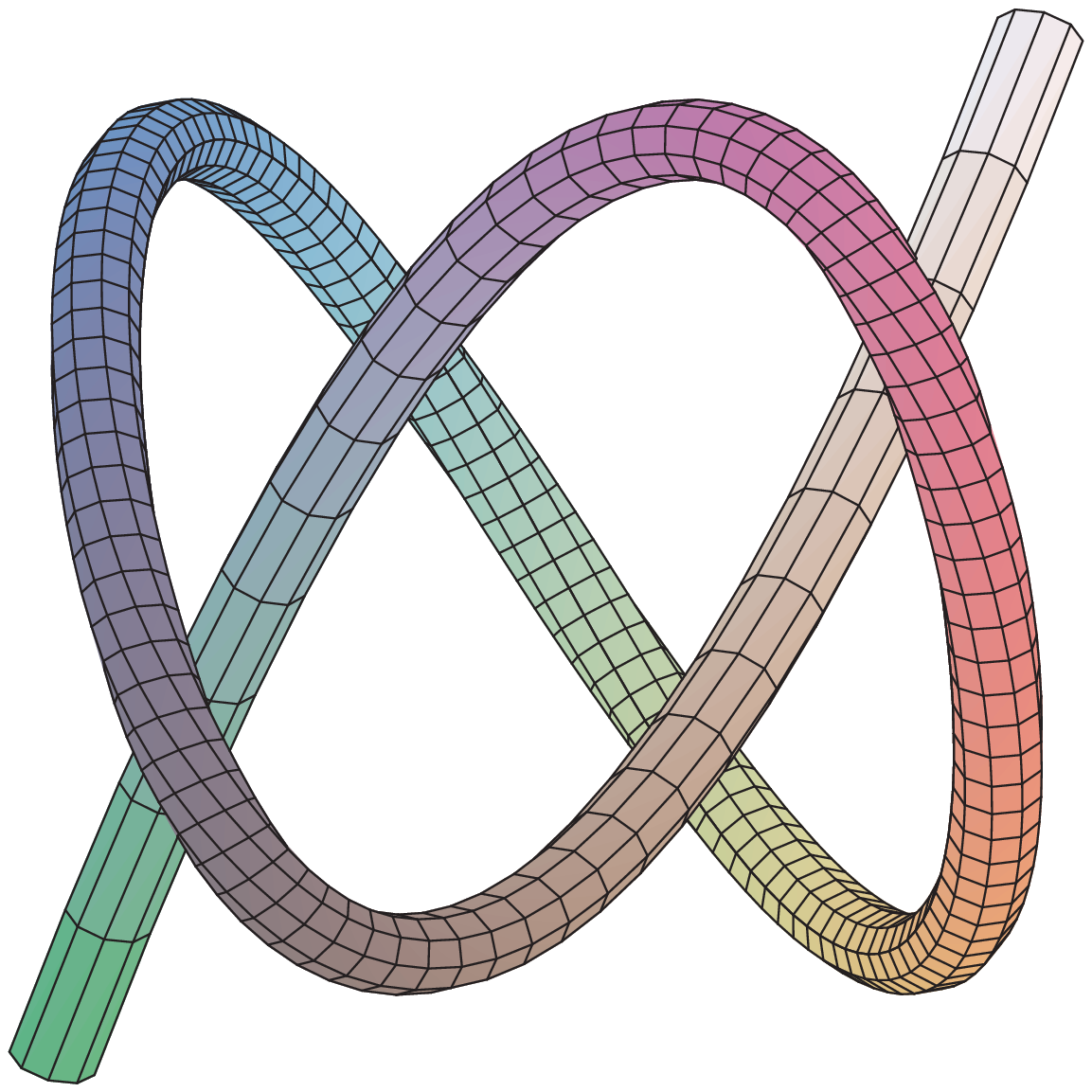,width=4cm}\\
bottom view of $K_3$&&face view of $K_3$
\end{tabular}
\end{center}
\vspace{-8pt}
\caption{The trefoil knot $K_3$\label{K3}}
\end{figure}
The plane curve $(T_3(t),T_4(t))$ has 3 crossing points. 
The plane curve $(T_3(t),T_5(t))$ has 4 crossing points corresponding to
parameters $(s_i,t_i)$ with 
$$s_1 < s_2 < s_3 <s_4 < t_2<t_1<t_4<t_3$$
so there do not exist real numbers   
$s_1<s_2<s_3$, and $t_1<t_2<t_3$ such that $ x(s_i)=x(t_i), \   z(s_i)=z(t_i).$

This example  shows that our  method cannot
be generalized when
the projections of $K_n$ have at least  $n+1$ crossing points.
\subsection{Parametrization of $K_5$}
Let us consider the curve of degree $(3,7,8)$: 
\begin{eqnarray*}
x &=&T_3(t), \\
y &=&T_8(t)-2\, T_6(t)+2.189\,T_4(t)- 2.170\, T_2(t),\\
z &=&T_7(t)-0.56\, T_5(t)- 0.01348 \,T_1(t).
\end{eqnarray*}
The curve $(x(t),y(t))$ has exactly 5 double points when the
projection $(x(t),z(t))$ has exactly 6. Note here that $\deg z(t) <
\deg y(t)$.
\begin{figure}[th]
\begin{center}
\begin{tabular}{ccc}
\psfig{file=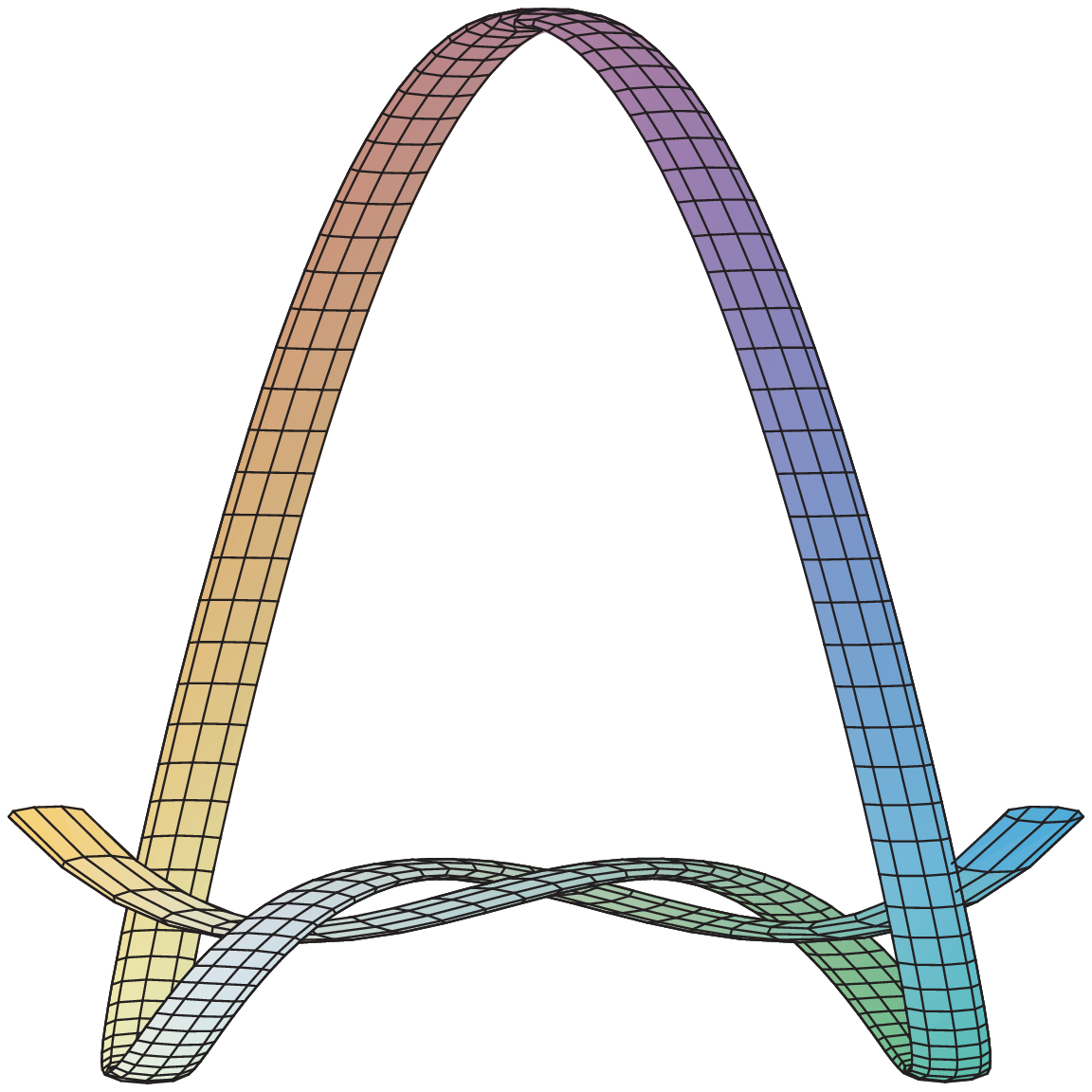,width=5cm}
&\hbox{}\hspace{1cm}\hbox{}
&\psfig{file=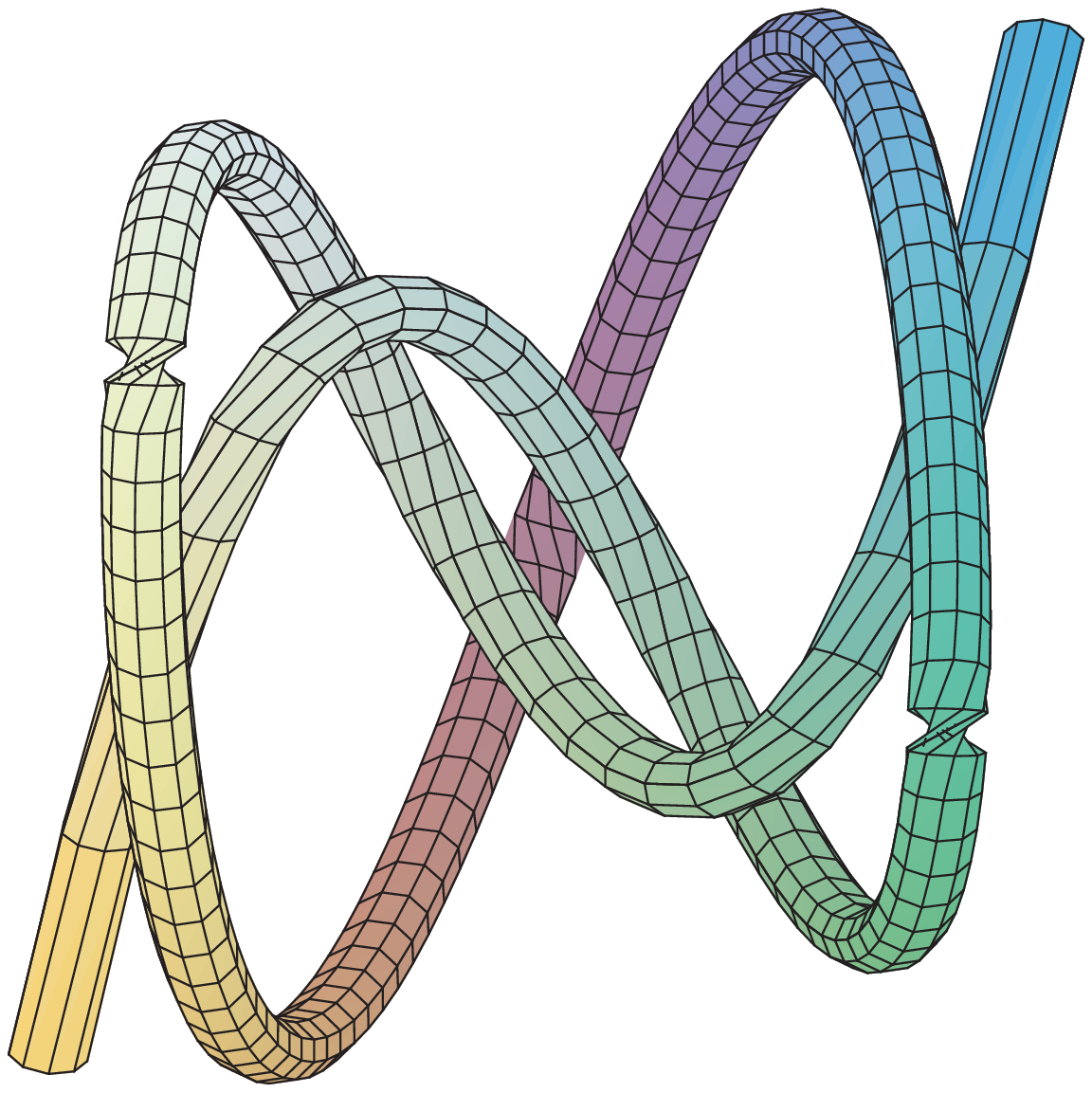,width=5cm}\\
bottom view&&face view
\end{tabular}
\end{center}
\vspace{-8pt}
\caption{Knot $K_5$\label{fig5}}
\end{figure}
In conclusion we have found a curve of degree $(3,7,8)$. Using our theorem,
we see that this curve has minimal degree. A. Ranjan and R. Mishra showed the
existence of such an example (\cite{RS,Mi}). 
\subsection{Parametrization of $K_7$}
We choose 
\begin{figure}[th]
\begin{center}
\begin{tabular}{ccc}
{\includegraphics[width=3cm]{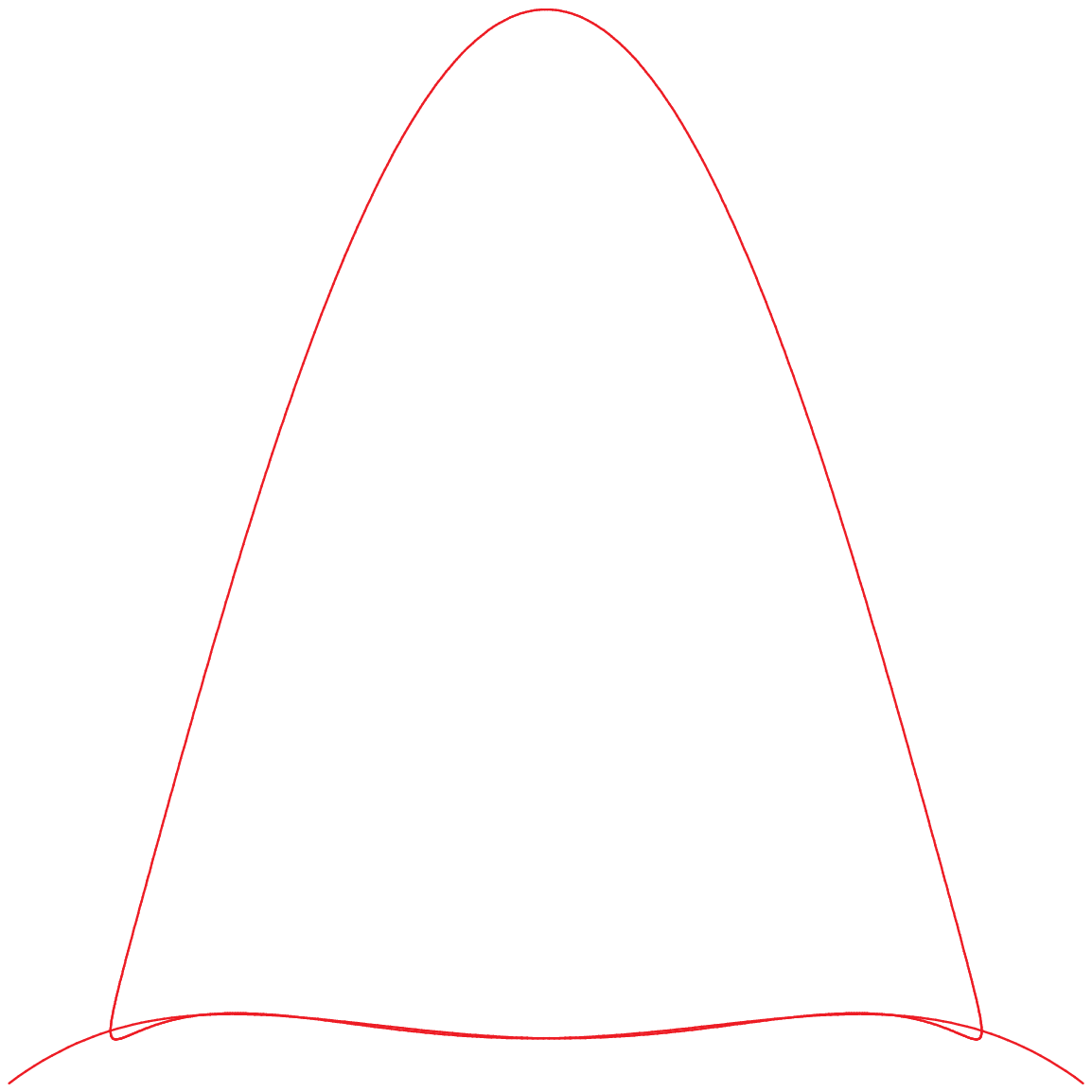}}&\quad&
{\includegraphics[width=7cm]{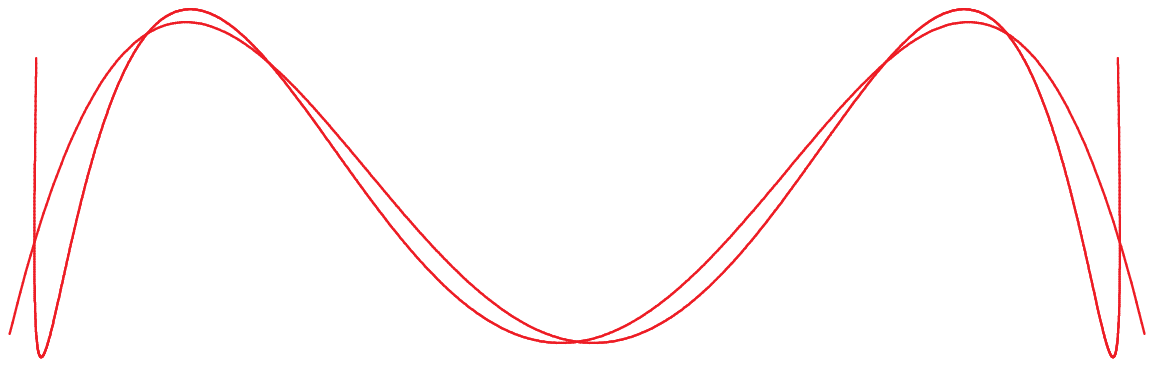}}\\
Bottom view&& Zoom on the bottom view
\end{tabular}
\end{center}
\vspace{-8pt}
\caption{Knot $K_7$\label{fig6}}
\end{figure}
\begin{eqnarray*}
x&=& T_{{3}} \left( t \right), \\
y&=& T_{{10}} \left( t \right) - 2.360\,T_{{8}} \left( t \right) + 4.108\,T_{{6}} \left( t \right) - 6.037\,T_{{4}}
 \left( t \right) + 7.397\,T_{{2}} \left( t \right),\\
z&=&T_{{11}}
 \left( t \right) + 3.580 \,T_{{7}} \left( t \right) - 3.739\,T_{{5}}
 \left( t \right) - T_{{1}} \left( t \right). 
\end{eqnarray*}
\parpic(0cm,7cm)(0cm,6cm)[lt]{\includegraphics[width=6cm]{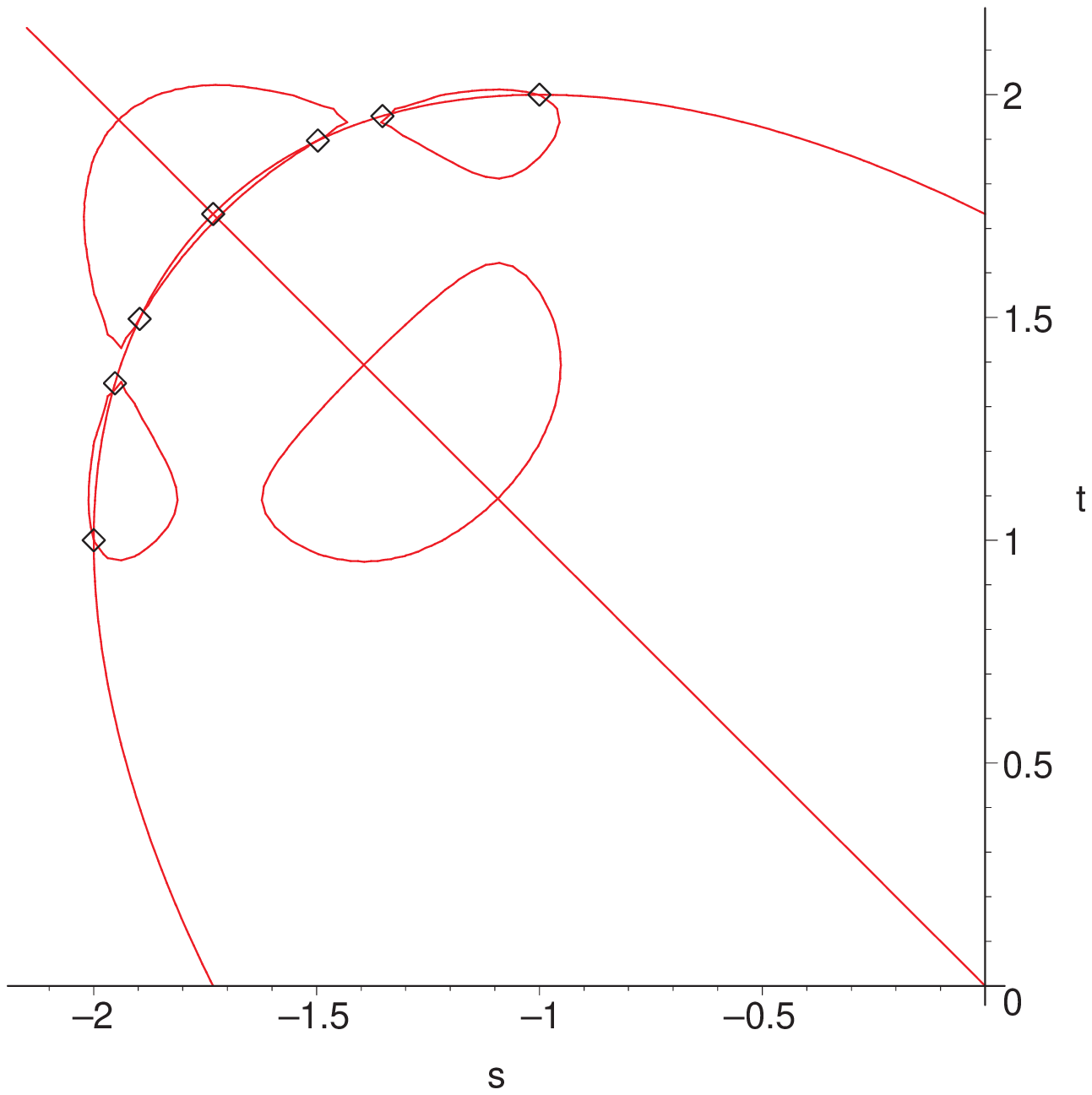}}
The values of the parameters corresponding to the double points are
  obtained as intersection points  between the ellipse $$t^2 + s^2 +st
  -3 =0$$ and the curve of degree 9: 
 $$(y(t)-y(s))/(t-s) = 0$$
The curve $(x(t),y(t))$ has exactly 7 double points 
corresponding to $\cos(\alpha)= \{ \pm 1/2, \pm 3/10,\pm 2/10,0\}$.

In conclusion we have found a curve of degree $(3,10,11)$. 
Using our theorem, this curve has minimal degree.

\subsection{Parametrization of $K_9$}
We choose polynomials of degree $(3,13,14)$.
\begin{eqnarray*}
x&=&T_{{3}} \left( t \right) ,\\
y&=& 
T_{{14}} \left( t \right) 
-  4.516\,T_{{12}} \left( t \right)  
+ 12.16\,T_{{10}} \left( t \right) 
- 24.46\,T_{{8}} \left( t \right) 
+ 39.92\,T_{{6}} \left( t \right) \\
& & 
\quad\quad- 55.30\,T_{{4}} \left( t \right)+66.60\,T_{{2}} \left( t \right) 
,\\
z&=&
T_{{13}} \left( t \right) 
- 2.389\,T_{{11}} \left( t \right)
- 5.161\,T_{{7}} \left( t \right) 
+ 5.161\,T_{{5}} \left( t \right) + 1.397\,T_{{1}} \left( t \right). 
\end{eqnarray*}
The curve $(x(t),y(t))$ has exactly 9 double points corresponding
to $\cos(\alpha)= \{ \pm 1/2, \pm 3/10,\pm 2/10, \pm 1/10,0\}$.
\begin{figure}[th]
\begin{center}
\begin{tabular}{ccc}
\psfig{file=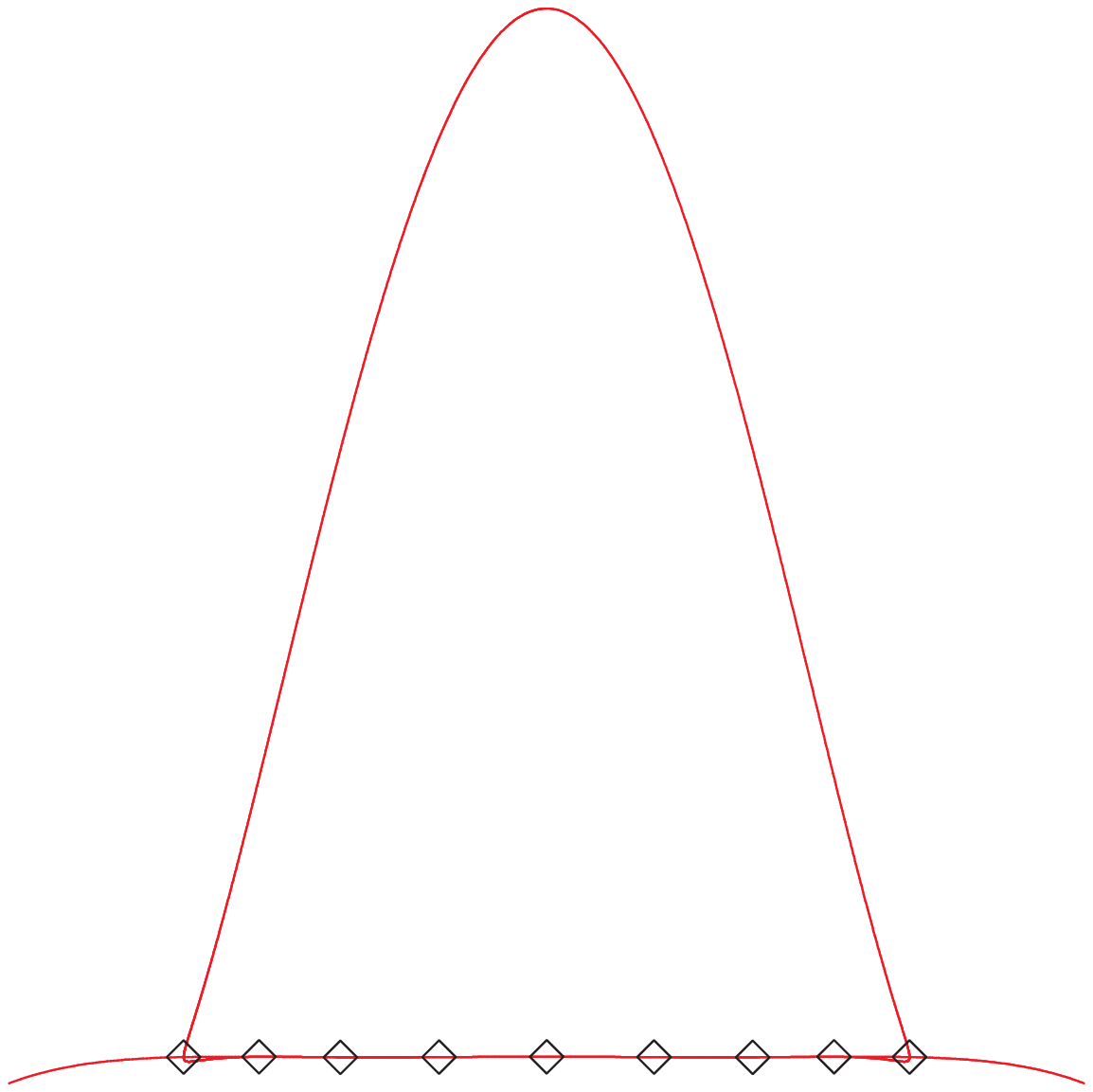,width=5cm}&&
\psfig{file=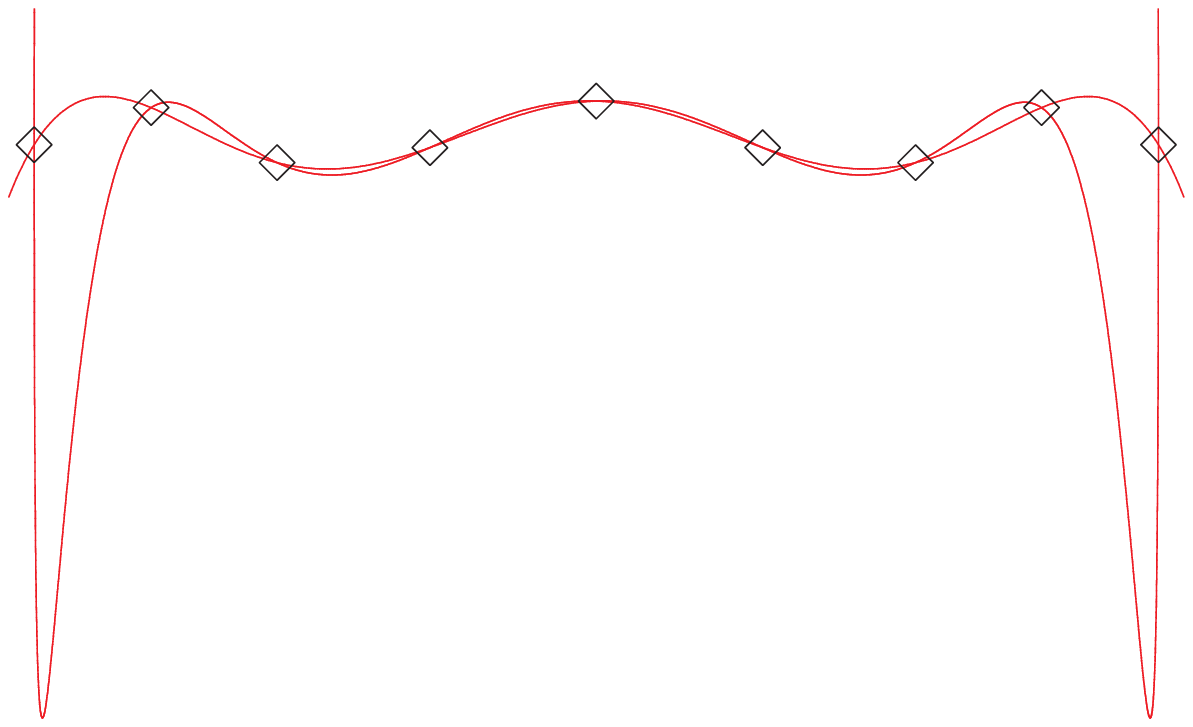,width=7cm}\\
Bottom view &  & Zoom on the bottom view
\end{tabular}
\end{center}
\vspace{-8pt}
\caption{Knot $K_9$\label{fig9}}
\end{figure}
One can prove that it is minimal under the assumption that the
projection $(x(t),y(t))$ has exactly 9 double points. 
\section*{Conclusion}
We have found minimal degree polynomial curves for torus knots $K_n$,
$n =3,5,7$. For degree $9$, one can prove that it is minimal under the
assumption that the projection $(x(t),y(t))$ has exactly 9 double
points. We have  similar constructions for higher degrees. 

\section*{Acknowledgments}
We would like to thank  Julien March{\'e} for fruitful discussions on
knot theory.

\end{document}